\documentclass[12pt]{article}
\usepackage{amsmath}
\usepackage{amssymb}
\usepackage{amscd}
\usepackage{amsthm}

\theoremstyle{plain}
\newtheorem{Thm}{Theorem}[section]

\newtheorem{Prop}[Thm]{Proposition}
\newtheorem{Cor}[Thm]{Corollary}
\newtheorem{Lem}[Thm]{Lemma}

\theoremstyle{definition}
\newtheorem{Defn}[Thm]{Definition}

\theoremstyle{Remark}

\numberwithin{equation}{section}

\title{On algebraic fiber spaces}
\author{Yujiro Kawamata}

\begin{document}

\maketitle

An algebraic fiber space is a relative version of an algebraic variety.
We prove some basic topological and analytical 
properties of algebraic fiber spaces.
We start with reviewing a result by Abramovich and Karu on the 
standard toroidal models of algebraic fiber spaces in \S1.
We can eliminate the singularities of the fibers in a topological sense
by performing the real oriented blowing-up on the toroidal model (\S2).
The rest of the paper concerns the Hodge theory of the degenerate fibers.
In \S3, we define the weight filtration on the cohomology with 
coefficients in $\mathbb{Z}$.
We prove the logarithmic and relative version of the Poincar\'e lemma 
in \S 4.
Finally we prove that the logarithmic de Rham complex gives rise to a 
cohomological mixed Hodge complex in \S 5 (Theorem~\ref{main}).  
As a corollary, we prove that certain spectral 
sequences degenerate.

We consider only varieties and morphisms defined over $\mathbb{C}$.
The topology is the classical (or Euclidean) topology instead of the 
Zariski topology unless stated otherwise.

%%%%%%%%%%%%%%%%%%%%%%%%%%%%%%%%%%%%%%%%%%%%%%%%%%%%%%%%%%%%%%%%%%%%%%%%

\section{Weak semistable model}

An {\it algebraic variety} in this paper is a reduced and irreducible scheme
of finite type over $\text{Spec }\mathbb{C}$.
An {\it algebraic fiber space} is a relative version of an algebraic variety;
it is a morphism $f: X \to Y$ of algebraic varieties which is generically 
surjective and such that the geometric generic fiber is 
reduced and irreducible.
We look for a standard model of an algebraic fiber space in the category 
of toroidal varieties.

\begin{Defn}\label{toroidal variety}
A {\it toric variety} $(V, D)$ is a pair consisting of a normal 
algebraic variety and a Zariski closed subset such that an algebraic 
torus $T \cong (\mathbb{C}^*)^n$ acts on $V$ 
with an open orbit $V \setminus D$.

A {\it toroidal variety} is a pair $(X, B)$ consisting of an algebraic 
variety and a Zariski closed subset which is locally analytically 
isomorphic to 
toric varieties in the following sense: for each point $x \in X$, 
there exists a toric variety $(V_x, D_x)$ with a fixed point $x'$, 
called a {\it local model} at $x$,
and open neighborhoods $U_x$ and $U'_{x'}$ of $x \in X$ and $x' \in V_x$ 
in the classical topology such that 
$(U_x, B \cap U_x)$ is isomorphic to $(U'_{x'}, D_x \cap U'_{x'})$.

A toroidal variety $(X, B)$ is called {\it strict} if any irreducible 
component of $B$ is normal.

A strict toroidal variety $(X, B)$ is called a {\it smooth toroidal variety} if
$X$ is smooth and $B$ is a simple normal crossing divisor.
Namely, a local model of a smooth toroidal variety has the form
$(\mathbb{C}^n, \text{div}(x_1 \dots x_{n'}))$, where $(x_1, \dots, x_n)$ are 
the coordinates.

A strict toroidal variety $(X, B)$ is called a {\it quasi-smooth toroidal 
variety} 
if its local model is a quotient of a smooth local model by a finite abelian
group action which is fixed point free on $X \setminus B$.
Namely, a local model has the form $(\mathbb{C}^n/G, 
\text{div}(x_1 \dots x_{n'})/G)$, where $G$ is a finite abelian
subgroup of $GL(n, \mathbb{C})$ 
which acts on $\mathbb{C}^n$ diagonally on the coordinates
$(x_1, \dots, x_{n'})$.

Let $(V, D)$ and $(W, E)$ be toric varieties with the actions of algebraic tori
$T$ and $S$.
A {\it toric morphism} $g: (V, D) \to (W, E)$ 
is a morphism $g: V \to W$ of algebraic varieties 
which is compatible with a homomorphism $g_0: T \to S$ of algebraic groups.

Let $(X, B)$ and $(Y, C)$ be toroidal varieties.
A {\it toroidal morphism} $f: (X, B) \to (Y, C)$ 
is a morphism $f: X \to Y$ of algebraic varieties 
such that $f(X \setminus B) \subset Y \setminus C$ and that
for any point $x \in X$ and any local model $(W_y, E_y)$ at $y = f(x) \in Y$,
there exists a local model $(V_x, D_x)$ at $x$ and a toric morphism 
$g: (V_x, D_x) \to (W_y, E_y)$ which is locally analytically isomorphic to $f$.
\end{Defn}

The resolution theorem of singularities by Hironaka implies 
that there is always a 
smooth birational model for any algebraic variety;
for any complete algebraic variety $X$,
there exists a birational morphism $\mu: Y \to X$ from a smooth 
projective variety.

As for the relative version of the resolution theorem, one cannot expect 
that we have a birational model which is a smooth morphism.
Indeed, we cannot eliminate singular fibers. 
Instead of the smooth model, Abramovich and Karu obtained 
a toroidal model by using the method of de Jong on the
moduli space of stable curves.

\begin{Thm}[Abramovich-Karu \cite{Abramovich-Karu}]\label{Abramovich-Karu}
Let $f_0: X_0 \to Y_0$ be a surjective morphism of complete algebraic 
varieties whose geometric generic fiber is irreducible.
Let $B_0$ and $C_0$ be Zariski closed subsets of $X_0$ and $Y_0$ such that
$f_0(X_0 \setminus B_0) \subset Y_0 \setminus C_0$.
Then there exist a quasi-smooth projective toroidal variety $(X, B)$,
a smooth projective toroidal variety $(Y, C)$, birational morphisms 
$\mu: X \to X_0$ and $\nu: Y \to Y_0$, and a toroidal and equi-dimensional 
morphism $f: (X, B) \to (Y, C)$ such that 
$\mu(X \setminus B) \subset X_0 \setminus B_0$, $\nu(Y \setminus C) \subset 
Y_0 \setminus C_0$ and $\nu \circ f = f_0 \circ \mu$. 
\end{Thm}

Let $f: (X, B) \to (Y, C)$ be a toroidal and equi-dimensional morphism of
quasi-smooth toroidal varieties.
We can describe $f$ explicitly by using local coordinates as follows.
Let us fix $x \in X$ and $y = f(x) \in Y$.
Let $n = \dim X$ and $m = \dim Y$.
We have local models of $X$ and $Y$: there are integers $0 \le n' \le n$ 
and $0 \le m' \le m$, finite abelian groups $G$ and $H$ which act diagonally
on the first $n'$ and $m'$ coordinates of 
the polydisks $\Delta^n = \{(x_1, \dots, x_n) \vert \vert x_i \vert < 1\}$
and $\Delta^m = \{(y_1, \dots, y_m) \vert \vert y_j \vert < 1\}$,
and open neighborhoods $U$ and $V$ of $x \in X$ and $y \in Y$ 
in the classical topology 
such that $(U, B \cap U) \cong (\Delta^n/G, \text{div}(x_1 \dots x_{n'})/G)$ 
and $(V, C \cap V) \cong (\Delta^m/H, \text{div}(y_1 \dots y_{m'})/H)$. 
We may assume that the fixed locus of each element of $G$ and $H$ except the
identities has codimension at least $2$.
Then the morphism $f$ induces a morphism 
$\tilde f: \Delta^n \to \Delta^m$.
The fact that $f$ is toroidal and equi-dimensional means the following:
there are integers $0 = t_0 < t_1 < t_2 < \dots < t_{m'} \le n'$ and 
$1 \le l_i$ ($i=1, \dots, t_{m'}$) such that
\begin{equation}\label{equation of weak semistable reduction}
\tilde f^*(y_j) = \prod_{k=1}^{t_j-t_{j-1}} x_{t_{j-1}+k}^{l_{t_{j-1}+k}} 
\end{equation}
for $j=1, \dots, m'$ and $\tilde f^*(y_j) = x_{n'-m'+j}$ for
$j=m'+1, \dots, m$.
Indeed, since $f$ is toroidal, it is expressed by monomials by some 
local coordinates.
The equi-dimensionality implies that the sets of indices $i$ of the $x_i$ 
on the right hand side for different $j$'s are disjoint. 

By using \cite{Characterization}, Abramovich and Karu obtained a model
with reduced fibers:

\begin{Cor}\label{weak semistable reduction}
Let $(X, B)$ be a quasi-smooth projective toroidal variety,
$(Y, C)$ a smooth projective toroidal variety, and 
$f: (X, B) \to (Y, C)$ a toroidal and equi-dimensional morphism
as in Theorem~\ref{Abramovich-Karu}.
Then there exists a finite and surjective morphism 
$\pi_Y: (Y', C') \to (Y, C)$ from a smooth projective toroidal variety 
such that $\pi^{-1}(C) = C'$ and the following holds: 
if we set $X'$ to be the normalization of the fiber product $X \times_Y Y'$ and
$B' = \pi_X^{-1}(B)$ for the induced morphism $\pi_X: X' \to X$,
then $(X', B')$ is a quasi-smooth projective toroidal variety 
and the induced morphism $f': (X', B') \to (Y', C')$ is a
toroidal and equi-dimensional morphism whose fibers are reduced.
\end{Cor}

The morphism $f'$ is called a {\it weak semistable reduction} of $f_0$.
We note that the fibers of $f'$ are reduced if and only if 
all the exponents $l_i$ are equal to $1$ in the corresponding local 
description (\ref{equation of weak semistable reduction}) for $f'$.

%%%%%%%%%%%%%%%%%%%%%%%%%%%%%%%%%%%%%%%%%%%%%%%%%%%%%%%%%%%%%%%%%%%%%%%%%

\section{Real oriented blowing-up}

It has been known that the general fiber of a semistable degeneration 
is topologically homeomorphic to the real oriented blowing-up of the 
singular fiber (e.g. \cite{Persson}).
\cite{Kawamata-Namikawa} used this knowledge 
to put a $\mathbb{Z}$-structure on a
certain cohomological mixed Hodge complex.
The real oriented blowing-up is a special case of the associated 
logarithmic topological space to a logarithmic complex space defined by 
\cite{Kato-Nakayama}.

\begin{Defn}\label{real oriented blowing-up}
The {\it real oriented blowing-up} of a quasi-smooth toroidal variety 
$(X, B)$ is a real analytic morphism from a 
real analytic manifold with boundary to a complex variety 
$\rho_X: X^{\#} \to X$ defined by the following recipe:

(0) If there is no boundary, then $\rho_X$ is the identity:
if $X = \Delta = \{z \in \mathbb{C} \vert \vert z \vert < 1\}$ and 
$B = \emptyset$,
then $X^{\#} = X$ and $\rho_X = \text{Id}$.

(1) This is the basic case: 
if $X = \Delta$ and 
$B = \{0\}$, then
$X^{\#} = [0,1) \times S^1$ and $\rho_X(r, \theta) = re^{i\theta}$.

(2) The real oriented blowing-up of a product is the 
product of the real oriented blowing-ups: 
if $X = X_1 \times X_2$ and $B = p_1^*B_1 + p_2^*B_2$, then
$X^{\#} = X_1^{\#} \times X_2^{\#}$ and
$\rho_X = \rho_{X_1} \times \rho_{X_2}$.

(3) The real oriented blowing-up of a quotient of a smooth toroidal variety
by a diagonal action of a finite abelian group which is free on the 
complement of the boundary divisor is the quotient of 
the real oriented blowing-up:
if $X = X_1/G$ and $B = B_1/G$, then
$X^{\#} = X_1^{\#}/G$ and $\rho_X = \rho_{X_1}/G$.
We note that the action of $G$ on $X_1^{\#}$ is free so that 
$X^{\#}$ has no singularities.

(4) We can glue together the real oriented blowing-ups of local models:
if $X = \bigcup_i X_i$ and $B = \bigcup_i B_i$, then
$X^{\#} = \bigcup_i X_i^{\#}$ and
$\rho_X = \bigcup_i \rho_{X_i}$.
\end{Defn}

\begin{Prop}\label{real oriented blowing-up and neighborhood}
Let $(X, B)$ be a quasi-smooth toroidal variety.
Then the real oriented blowing-up $X^{\#}$ is homeomorphic to
the complement of an $\epsilon$ neighborhood of the boundary $B$ in $X$
with respect to some metric for sufficiently small $\epsilon$.
\end{Prop}

The real oriented blowing-up is functorial:

\begin{Prop}\label{real oriented blowing-up is functorial}
Let $f: (X, B) \to (Y, C)$ be a toroidal morphism between
quasi-smooth toroidal varieties.
Let $\rho_X: X^{\#} \to X$ and $\rho_Y: Y^{\#} \to Y$
be the real oriented blowing-ups.
Then a morphism of real analytic varieties with boundaries
$f^{\#}: X^{\#} \to Y^{\#}$ is induced so that
the following diagram is commutative:
\[
\begin{CD}
X^{\#} @>{\rho^1_X}>> X \times_Y Y^{\#} @>{\text{pr}_1}>> X \\
@V{f^{\#}}VV            @V{\text{pr}_2}VV    @VfVV \\
Y^{\#} @>=>>            Y^{\#}            @>{\rho_Y}>> Y.
\end{CD}
\]
\end{Prop}

By using the real oriented blowing-up, we can eliminate the singularities of 
fibers topologically:

\begin{Thm}\label{locally topologically trivial}
Let $f: (X, B) \to (Y, C)$ be a proper surjective toroidal and equi-dimensional
morphism of quasi-smooth toroidal varieties.
Then the induced morphism 
$f^{\#}: X^{\#} \to Y^{\#}$
is locally topologically trivial in the following sense:
each point $y' \in Y^{\#}$ has an open neighborhood $V'$ such that 
$(f^{\#})^{-1}(V')$ is homeomorphic to $V' \times (f^{\#})^{-1}(y')$ 
over $V'$.
\end{Thm}

\begin{proof}
We use the local description explained after Theorem~\ref{Abramovich-Karu}.
We write $\rho_{\tilde X}^*(x_i) = r_ie^{\sqrt{-1}\theta_i}$ 
($1 \le i \le n'$) and 
$\rho_{\tilde Y}^*(y_j) = s_je^{\sqrt{-1}\phi_j}$ ($1 \le i \le m'$).
Then the map $f^{\#}$ is described by the following formulas: 
\begin{equation}\label{equation of weak semistable reduction2}
\begin{split}
(\tilde f^{\#})^*s_j &= \prod_{k=1}^{t_j-t_{j-1}} 
r_{t_{j-1}+k}^{l_{t_{j-1}+k}} \\
(\tilde f^{\#})^*\phi_j &= \sum_{k=1}^{t_j-t_{j-1}} 
l_{t_{j-1}+k}\theta_{t_{j-1}+k}
\end{split}
\end{equation}
for $j=1, \dots, m'$ and $(\tilde f^{\#})^*(y_j) = x_{n'-m'+j}$ for
$j=m'+1, \dots, m$.
We note that the actions of $G$ and $H$ are rivial on the $r_i$ and $s_j$ while
those on the $\theta_i$ and $\phi_j$ are translations.
Since 
\[
\{(r_1, \dots, r_t) \vert t_i \in [0,1), \prod_i t_i = c\} \cong [0,1)^{t-1}
\]
for any $c \in [0,1)$, 
the fibers of $f^{\#}: U^{\#} \to V^{\#}$ are 
homeomorphic to 
\[
[0,1)^{n'-m'} \times (S^1)^{n'-m'} \times (D^2)^{n-n'-m+m'}.
\]
The restriction of $f^{\#}$ on $U^{\#} \cap (f^{\#})^{-1}(V')$ 
for sufficiently small $V'$ 
is homeomorphic to the first projection of $V' \times (f^{\#})^{-1}(y')$ 
to $V'$.
By gluing together, we obtain our result.
\end{proof}

\begin{Cor}\label{locally constant sheaf}
$R^p(f^{\#})_*\mathbb{Z}_{X^{\#}}$ is a locally
constant sheaf on $Y^{\#}$ for any $p \ge 0$.
\end{Cor}

\begin{Prop}\label{real oriented blowing-up of a stratum}
Let $(X, B)$ be a quasi-smooth toroidal variety and let
$B = \sum_{i=1}^N B_i$ be the irreducible decomposition.
Let $E$ be a connected component of the intersection $\bigcap_{i=1}^l B_i$
of codimension $l$, 
and let $G = \sum_{i=l+1}^N G_i$ for $G_i = B_i \cap E$.
Then $(E, G)$ is a quasi-smooth toroidal variety. 
Let $\rho_X: X^{\#} \to X$ and 
$\rho_E: E^{\#} \to E$
be the real oriented blowing-ups.
Then $\rho_X$ induces a map
$\rho_X': \rho_X^{-1}(E) \to E^{\#}$ which is an $l$-times fiber 
product of oriented $S^1$-bundles.
\end{Prop}

\begin{proof}
Since the boundary $B$ has no self-intersection, $E$ is normal, hence
$(E, G)$ is toroidal and quasi-smooth.
$E \setminus G$ is smooth and the normal bundle
$N_{(E \setminus G)/(X \setminus G)}$ is the direct product of line bundles.
Since $E^{\#}$ is homeomorphic to the complement of an 
$\epsilon$ neighborhood of $G$ in $E$, we obtain our assertion.
\end{proof}

\begin{Cor}\label{direct image over Z}
\[
\begin{split}
&R^p(\rho_X)_*\mathbb{Z}_{X^{\#}}
\cong \bigwedge^p (\bigoplus_{i=1}^N \mathbb{Z}_{B_i}) \\
&R^p(\rho_X)_*\mathbb{Z}_{\rho_X^{-1}(E)}
\cong \bigwedge^p (\mathbb{Z}_E^l \oplus
\bigoplus_{i=l+1}^N \mathbb{Z}_{G_i})
\end{split}
\]
where the exterior products are taken respectively as 
$\mathbb{Z}_X$-modules and $\mathbb{Z}_E$-modules.
\end{Cor}

In the case $p = 0$, the above formula means that
$(\rho_X)_*\mathbb{Z}_{X^{\#}} \cong \mathbb{Z}_X$ and 
$(\rho_X)_*\mathbb{Z}_{\rho_X^{-1}(E)} \cong \mathbb{Z}_E$.
We note that $G_i$ may be empty or reducible.

%%%%%%%%%%%%%%%%%%%%%%%%%%%%%%%%%%%%%%%%%%%%%%%%%%%%%%%%%%%%%%%%%%%%%%%%%%%%%

\section{Weight filtration}

We put a weight filtration on a complex on singular fibers.
The definition of the filtration is natural 
thanks to the geometric construction of the real oriented blowing-up.
One can compare with rather complicated earlier 
definitions in \cite{Steenbrink}
and \cite{Fujisawa}.  

In this section, we denote by $f: (X, B) \to (Y, C)$ 
a weak semistable reduction as in Corollary~\ref{weak semistable reduction}.
We fix $y \in Y$ and let $E_1, \dots, E_l$ be irreducible components
of the fiber $E = f^{-1}(y)$.
We have $\dim E_k = n - m$ for any $k$.
For any combination of integers $1 \le i_0 < \dots < i_t \le l$, 
we define $E_{i_0, \dots, i_t} = \bigcap_{k=0}^t E_{i_k}$.
We note that $\dim E_{i_0, \dots, i_t}$ may be larger or smaller than
$n - m - t$.
For example, if a local model $f: \Delta^4 \to \Delta^2$ is given by
$f^*(y_1) = x_1x_2$ and $f^*(y_2) = x_3x_4$ with $y = (0,0)$, then
the irreducible components of $E$ are 
$E_1 = \{x_1=x_3=0\}$, $E_2 = \{x_1=x_4=0\}$, $E_3 = \{x_2=x_3=0\}$
and $E_4 = \{x_2=x_4=0\}$ so that
$E_{14} = E_{23} = E_{234} = E_{134} = E_{124} = E_{123} = E_{1234}$.  

We call each $E_{i_0, \dots, i_t}$ a {\it stratum} of $E$. 
Let $G_{i_0, \dots, i_t}$ be the union of all the strata which are
properly contained in $E_{i_0, \dots, i_t}$.
Then the pair 
\[
(E_{i_0, \dots, i_t}, G_{i_0, \dots, i_t})
\]
is again a quasi-smooth toroidal variety.
Let $\rho_{i_0, \dots, i_t}: E_{i_0, \dots, i_t}^{\#} \to 
E_{i_0, \dots, i_t}$ be the real oriented blowing-up.

Let $E^{[t]}$ be the disjoint union of all the $E_{i_0, \dots, i_t}$.
Then we have an exact sequence
\[
0 \to \mathbb{Z}_E \to \mathbb{Z}_{E^{[0]}} \to
\mathbb{Z}_{E^{[1]}} \to \dots \to \mathbb{Z}_{E^{[t]}} \to \dots 
\]
where we identified the sheaves $\mathbb{Z}_{E^{[t]}}$
with their direct image sheaves on $E$ and the differentials are 
given by the alternate sums of the restriction maps.

If the number of irreducible components of $C$ which contain $y$ is $m'$,
then $\rho_Y^{-1}(y)$ is homeomorphic to $(S^1)^{m'}$.
Let $D = \rho_X^{-1}(E) = (f^{\#})^{-1}\rho_Y^{-1}(y)$. 
We have natural maps $\rho_X: D \to E$ and 
$f^{\#}: D \to \rho_Y^{-1}(y)$.
Corresponding to the irreducible decomposition of $E$, we have
$D = \bigcup_{i=1}^l D_i$ for $D_i = \rho_X^{-1}(E_i)$.
We write $D_{i_0, \dots, i_t} = \bigcap_{k=0}^t D_{i_k}$ and
$D^{[t]} = \bigsqcup_{1 \le i_0 < \dots < i_t \le l} D_{i_0, \dots, i_t}$.
Then we have similarly an exact sequence
\[
0 \to \mathbb{Z}_D \to \mathbb{Z}_{D^{[0]}} \to \mathbb{Z}_{D^{[1]}}
\to \dots \to \mathbb{Z}_{D^{[t]}} \to \dots. 
\]
We fix a point $\bar y \in \rho_Y^{-1}(y)$ and let 
$\bar D = (f^{\#})^{-1}(\bar y)$.
We set $\bar D_{i_0, \dots, i_t} = \bar D \cap D_{i_0, \dots, i_t}$
and $\bar D^{[t]} = 
\bigsqcup_{1 \le i_0 < \dots < i_t \le l} \bar D_{i_0, \dots, i_t}$.

\begin{Prop}\label{direct product}
Let $t' = \dim E - \dim E_{i_0, \dots, i_t}$.
Then $D_{i_0, \dots, i_t}$ is 
homeomorphic to the direct product 
$\bar D_{i_0, \dots, i_t} \times \rho_Y^{-1}(y)$, and $f^{\#}$ corresponds
to the second projection.
Moreover, $\bar D_{i_0, \dots, i_t}$ is a $t'$-times
fiber product of $S^1$-bundles over $E_{i_0, \dots, i_t}^{\#}$.
\end{Prop}

\begin{proof}
We assume first that $t = t' = 0$.
Then we have $\bar D_{i_0} = E_{i_0}^{\#}$.
By Proposition~\ref{real oriented blowing-up of a stratum}, 
$\rho_X$ induces an $(S^1)^{m'}$-bundle 
$\rho'_X: D_{i_0} \to E_{i_0}^{\#}$, and
we have a map
\[
\rho_X' \times f^{\#}: D_{i_0} \to 
E_{i_0}^{\#} \times \rho_Y^{-1}(y).
\]
Since the fiber of $f$ is reduced, we may take the homeomorphism to the 
$\epsilon$ neighborhood as in Proposition~\ref{real oriented blowing-up of a 
stratum} in a suitable way to conclude that $\rho_X' \times f^{\#}$
is bijective and a homeomorphism.

In the general case, if we restrict $\rho_X' \times f^{\#}$ to 
$D_{i_0, \dots, i_t}$, we 
obtain a homeomorphism 
\[
D_{i_0, \dots, i_t} \to 
\rho_{i_0}^{-1}(E_{i_0, \dots, i_t}) \times \rho_Y^{-1}(y).
\]
By applying Proposition~\ref{real oriented blowing-up of a 
stratum} again to $E_{i_0, \dots, i_t} \subset E_{i_0}$, 
we conclude the proof.
\end{proof}

\begin{Defn}\label{weight filtration}
We define the {\it weight filtration} $W_q$ on 
the complex of sheaves $R(\rho_X^1)_*\mathbb{Z}_D$
on $E \times \rho_Y^{-1}(y)$ by
\[
\begin{split}
&W_q(R(\rho_X^1)_*\mathbb{Z}_D) \\
&= (\tau_{\le q}(R(\rho_X^1)_*\mathbb{Z}_{D^{[0]}}) \to 
\tau_{\le q+1}(R(\rho_X^1)_*\mathbb{Z}_{D^{[1]}}) \to \\
&\quad \dots \to \tau_{\le q+t}(R(\rho_X^1)_*\mathbb{Z}_{D^{[t]}}) \to \dots)
\end{split}
\]
in the derived category of complexes of sheaves on 
$E \times \rho_Y^{-1}(y)$, where $\tau$ denotes the truncation
(cf. 1.4.6 of \cite{Deligne}).
\end{Defn}

If we take the successive quotients with respect to this filtration, 
then we can check that all the differentials of the resolution become trivial:

\begin{Prop}\label{Gr}
\[
\begin{split}
&\text{Gr}_q^W(R(\rho_X^1)_*\mathbb{Z}_D) \\
&\cong R^q(\rho_X^1)_*\mathbb{Z}_{D^{[0]}}[-q]
\oplus R^{q+1}(\rho_X^1)_*\mathbb{Z}_{D^{[1]}}[-q-2] \oplus \\
&\quad \dots \oplus R^{q+t}(\rho_X^1)_*\mathbb{Z}_{D^{[t]}}[-q-2t]
\oplus \dots 
\end{split}
\]
\end{Prop}

\begin{Cor}
The local monodromies of $R^p(f^{\#})_*\mathbb{Z}_{X^{\#}}$ 
on $Y^{\#}$ are unipotent for any $p \ge 0$.
\end{Cor}

\begin{proof}
We consider the restriction of the locally constant sheaf
$R^p(f^{\#})_*\mathbb{Z}_{X^{\#}}$ on $\rho_Y^{-1}(y)$.
Since $f^{\#} = \text{pr}_2 \circ \rho_X^1$, 
the weight filtration induces a spectral sequence
\[
E_1^{p,q} = R^{p+q}(\text{pr}_2)_*(\text{Gr}_{-p}^W(R(\rho_X^1)_*\mathbb{Z}_D))
\Rightarrow R^{p+q}(f^{\#})_*\mathbb{Z}_{X^{\#}}.
\]
By Proposition~\ref{direct product}, 
$R^p(\rho_X^1)_*\mathbb{Z}_{D^{[t]}}$ 
is a constant sheaf for any $p,t$.
Hence we have our result.
\end{proof}

%%%%%%%%%%%%%%%%%%%%%%%%%%%%%%%%%%%%%%%%%%%%%%%%%%%%%%%%%%%%%%%%%%%%%%%%%%

\section{Relative log de Rham complex}

Let $(X, B)$ be a quasi-smooth toroidal variety of dimension $n$.
Then the sheaf of logarithmic $1$-forms $\Omega^1_X(\log B)$ is a locally
free sheaf of rank $n$.
Indeed, if the local model of $(X, B)$ at $x \in X$ is the
quotient of type $(\Delta^n/G, B_{n'}/G)$ with 
$B_{n'} = \text{div}(x_1 \dots x_{n'})$,
then the action of $G$ on the basis $dx_i/x_i$ ($1 \le i \le n'$) and 
$dx_i$ ($n' < i \le n$)
of the sheaf $\Omega^1_{\Delta^n}(\log B_{n'})$ is trivial.
The sheaf of logarithmic $p$-forms $\Omega^p_X(\log B) = 
\bigwedge^p \Omega^1_X(\log B)$ 
is a locally free sheaf of rank $\binom{n}{p}$.

Let $f: (X, B) \to (Y, C)$ be a toroidal and equi-dimensional morphism of
quasi-smooth toroidal varieties.
Let $\dim X = n$ and $\dim Y = m$.
Then
\[
\Omega^1_{X/Y}(\log) = \Omega^1_X(\log B)/f^*\Omega^1_Y(\log C)
\]
is a locally free sheaf of rank $n - m$.

Following \cite{Steenbrink} and \cite{Kato-Nakayama}, 
we define the structure sheaf $\mathcal{O}_{X^{\#}}$ of the
real oriented blowing-up by adding the logarithms of the coordinates:

\begin{Defn}\label{structure sheaf}
(0) If there is no boundary, then $\mathcal{O}_{X^{\#}}
= \mathcal{O}_X$.

(1) If $(X, B) = (\Delta, \{0\})$, then
\[
\mathcal{O}_{X^{\#}}
= \sum_{k \in \mathbb{Z}_{\ge 0}} 
\rho_X^{-1}(\mathcal{O}_X)(\log x)^k
\] 
where $\rho_X^{-1}$ denotes the inverse image for the sheaves of 
abelian groups, and $x$ is the coordinate.
The symbol $\log x$ is identified with a local section of 
$\rho_X^{-1}(\mathcal{O}_X)$ on $X \setminus B$.
We note that the multivalued function $\log x$ on $X$ becomes locally 
single valued on $X^{\#}$.
Therefore, the stalk of $\mathcal{O}_{X^{\#}}$ at $(r, \theta)$ is 
isomorphic to $\mathcal{O}_{X, x}$ if $r \ne 0$ and
to $\mathcal{O}_{X, 0} \otimes \mathbb{Z}[t]$ if $r = 0$, where $t$ is an 
independent variable corresponding to $\log x$.
We note that a section of $\mathcal{O}_{X^{\#}}$ is a finite polynomial 
on $\log x$ instead of an infinite power series.

(2) If $(X, B) = (X_1 \times X_2, p_1^*B_1 + p_2^*B_2)$, then
\[
\begin{split}
\mathcal{O}_{X^{\#}} = 
&(\text{pr}_1^{-1}(\mathcal{O}_{X_1^{\#}}) 
\otimes_{\text{pr}_1^{-1}\rho_{X_1}^{-1}(\mathcal{O}_{X_1})}
\rho_X^{-1}(\mathcal{O}_X)) \\
&\otimes_{\rho_X^{-1}(\mathcal{O}_X)}
(\text{pr}_2^{-1}(\mathcal{O}_{X_2^{\#}}) 
\otimes_{\text{pr}_1^{-1}\rho_{X_1}^{-1}(\mathcal{O}_{X_2})}
\rho_X^{-1}(\mathcal{O}_X)).
\end{split}
\]
For example, if $X = \Delta^n$ with $B = \text{div}(x_1 \dots x_{n'})$, 
then we have
\[
\mathcal{O}_{X^{\#}}
= \sum_{k_1, \dots, k_{n'} \in \mathbb{Z}_{\ge 0}}
\rho_X^{-1}(\mathcal{O}_X) \prod_{i=1}^{n'} (\log x_i)^{k_i}.
\]
In this case, the stalk of $\mathcal{O}_{X^{\#}}$ at any point over 
the origin of $X$ 
is isomorphic to $\mathcal{O}_{X, 0} \otimes \mathbb{Z}[t_1, \dots, t_{n'}]$,
where the $t_i$ are independent variables corresponding to the $\log x_i$.

(3) If $(X, B) = (X_1/G, B_1/G)$, then
$\mathcal{O}_{X^{\#}} = (\mathcal{O}_{X_1^{\#}})^G$.

(4) We can glue together the structure sheaves of local models.
\end{Defn}

\begin{Lem}\label{locally constant sheaf $L$}
Let $(X, B)$ be a quasi-smooth toroidal variety, $x \in X$ a point, and
let $n'$ be the number of irreducible components $B_i$ of $B$ which 
contain $x$.
Then $\mathcal{O}_{X^{\#}} \otimes_{\rho_X^{-1}(\mathcal{O}_X)}
\rho_X^{-1}(\mathbb{C}_x)$
is a locally constant sheaf $L$ on $\rho_X^{-1}(x) \cong (S^1)^{n'}$
whose fibers are isomorphic to $\mathbb{C}[t_1, \dots, t_{m'}]$,
where the $t_j$ are independent variables corresponding to the logarithms
of the local equations of the $B_i$ at $x$.
The monodromies $M_i$ of $L$ around the loop in $\rho_X^{-1}(x)$ 
corresponding to the $B_i$ 
are given by the following formula
\[
M_i(t_j) = t_j + 2\delta_{ij}\pi\sqrt{-1}.
\]
\end{Lem}

\begin{Defn}\label{log de Rham complex}
We define the {\em log de Rham complex} 
on the real blowing-up by
\[
\Omega^p_{X^{\#}} =
\rho_X^{-1}(\Omega^p_X(\log B)) \otimes_{\rho_X^{-1}(\mathcal{O}_X)} 
\mathcal{O}_{X^{\#}}
\]
for $p \ge 0$.
The differential of the complex $\Omega^{\bullet}_{X^{\#}}$
is defined by the rule: $d\log x = dx/x$. 
The {\it relative log de Rham complex} is defined by
\[
\Omega^p_{X^{\#}/Y^{\#}} =
\rho_X^{-1}(\Omega^p_{X/Y}(\log)) \otimes_{\rho_X^{-1}(\mathcal{O}_X)}
\mathcal{O}_{X^{\#}}.
\]
The differential of the complex 
$\Omega^{\bullet}_{X^{\#}/Y^{\#}}$ 
is induced from that of $\Omega^{\bullet}_{X^{\#}}$.
\end{Defn}

If $B = \emptyset$, then the Poincar\'e lemma says that
\[
\begin{split}
\mathbb{C}_X &\cong \Omega^{\bullet}_X \\
&= (\mathcal{O}_X \to \Omega^1_X \to \dots \to \Omega^p_X \to \dots)
\end{split}
\]
in the derived category of sheaves on $X$. 
We have the following logarithmic version:

\begin{Thm}\label{log Poincare}
(1) 
\[
\mathbb{C}_{X^{\#}} \cong \Omega^{\bullet}_{X^{\#}}
\]
in the derived category of sheaves on $X^{\#}$.

(2)
\[
R(\rho_X)_*\Omega^p_{X^{\#}} \cong \Omega^p_X(\log B)
\]
in the derived category of sheaves on $X$.

(3)
\[
R(\rho_X)_*\mathbb{C}_{X^{\#}} 
\cong \Omega^{\bullet}_X(\log B) 
\]
in the derived category of sheaves on $X$.
\end{Thm}

\begin{proof}
(1) This is a special case of Theorem 3.6 of \cite{Kato-Nakayama}.
For example, we check the case where
$(X, B) = (\Delta, \{0\})$.
Let $h = \sum_k h_k(\log x)^k$ for $h_k \in \mathcal{O}_{X,0}$.
If 
\[
dh = dx/x \sum_k (xh'_k+(k+1)h_{k+1})(\log x)^k = 0
\]
then we have $xh'_k+(k+1)h_{k+1} = 0$ for all $k$.
Since $h_k = 0$ for $k \gg 0$, we conclude that 
$h \in \mathbb{C}$ by the descending induction on $k$.

On the other hand, we can solve the equations 
\[
xg'_k + (k+1)g_{k+1} = h_k
\]
by the descending induction on $k$
to find $g = \sum_k g_k(\log x)^k$ such that $hdx/x = dg$.

(2) We check the case where $(X, B) = (\Delta, \{0\})$ and $p = 0$, 
that is, $(\rho_X)_*\mathcal{O}_{X^{\#}} 
= \mathcal{O}_X$ and $R^1(\rho_X)_*\mathcal{O}_{X^{\#}} = 0$.  
The case $p > 0$ follows from the projection formula,
and the general case is similarly proved.

The first equality follows from the fact that $\log x$ is multi-valued on
$X$ near $0$.
We have an exact sequence
\[
0 \to \rho_X^{-1}(\mathcal{O}_X) \to \mathcal{O}_{X^{\#}} \to
\bar L \otimes \mathcal{O}_{X,0} \to 0
\]
where $\bar L = L/\mathbb{C}$ is a locally constant 
sheaf on $\rho_X^{-1}(0) \cong S^1$ 
such that the stalk $\bar L_0$ at $0 \in S^1$ 
has a $\mathbb{C}$-basis $\{e_k\}_{k \in \mathbb{Z}_{>0}}$ 
corresponding to the $(\log x)^k$ with the monodromy action given by 
\[
M(e_k) = \sum_{j=1}^k \binom{k}{j}(2\pi\sqrt{-1})^{k-j}e_j.
\]
Then $H^0(S^1, \bar L) = \mathbb{C}e_1$ and $H^1(S^1, \bar L) = 0$. 
Hence 
\[
\begin{split}
&(\rho_X)_*\rho_X^{-1}(\mathcal{O}_X) \cong 
(\rho_X)_*\mathcal{O}_{X^{\#}} \\
&H^0(S^1, \bar L) \otimes \mathcal{O}_{X, 0}
\cong R^1(\rho_X)_*\rho_X^{-1}(\mathcal{O}_X) 
\end{split}
\]
and $R^1(\rho_X)_*\mathcal{O}_{X^{\#}} = 0$. 

(3) is a consequence of (1) and (2).
This follows also from Proposition 3.1.8 of \cite{Deligne}, because 
we have 
\[
R(\rho_X)_*\mathbb{C}_{X^{\#}}
\cong Rj_*\mathbb{C}_{X \setminus B}
\]
for the inclusion $j: X \setminus B \to X$
by Proposition~\ref{real oriented blowing-up and neighborhood}.
\end{proof}

In the relative setting, if $B = C = \emptyset$, then we have
$f^{-1}(\mathcal{O}_Y) \cong \Omega^{\bullet}_{X/Y}$ 
in the derived category of sheaves on $X$.
The logarithmic version is the following:

\begin{Thm}\label{relative log Poincare}
(1) 
\[
(f^{\#})^{-1}(\mathcal{O}_{Y^{\#}})
\cong \Omega^{\bullet}_{X^{\#}/Y^{\#}}
\]
in the derived category of sheaves on $X^{\#}$.

(2)
\[
R(\rho_X^1)_*\Omega^p_{X^{\#}/Y^{\#}}
\cong (\text{pr}_1)^{-1}(\Omega^p_{X/Y}(\log)) 
\otimes_{(\text{pr}_2)^{-1}\rho_Y^{-1}(\mathcal{O}_Y)}
(\text{pr}_2)^{-1}(\mathcal{O}_{Y^{\#}})
\]
in the derived category of sheaves on $X \times_Y Y^{\#}$.

(3) 
\[
R(\rho_X^1)_*(f^{\#})^{-1}(\mathcal{O}_{Y^{\#}})
\cong (\text{pr}_1)^{-1}(\Omega^{\bullet}_{X/Y}(\log)) 
\otimes_{(\text{pr}_2)^{-1}\rho_Y^{-1}(\mathcal{O}_Y)}
(\text{pr}_2)^{-1}(\mathcal{O}_{Y^{\#}})
\]
in the derived category of sheaves on $X \times_Y Y^{\#}$.
\end{Thm}

\begin{proof}
(1) We check the case where $(X, B) = (\Delta^2, \text{div}(x_1x_2))$, 
$(Y, C) = (\Delta, \{0\})$, and $f^*(y) = x_1x_2$.
The general case is similar.
Let 
\[
h = \sum_{k_1,k_2} h_{k_1,k_2}(\log x_1)^{k_1}(\log x_2)^{k_2}
\] 
for $h_{k_1,k_2} \in \mathcal{O}_{X,0}$.
If 
\[
\begin{split}
dh = &dx_1/x_1 \sum_{k_1,k_2} (x_1h_{k_1,k_2,x_1} - x_2h_{k_1,k_2,x_2} \\
&+(k_1+1)h_{k_1+1,k_2}-(k_2+1)h_{k_1,k_2+1})(\log x_1)^{k_1}(\log x_2)^{k_2}
= 0
\end{split}
\]
in $\Omega^1_{X^{\#}/Y^{\#}}$, then we have
\[
x_1h_{k_1,k_2,x_1} - x_2h_{k_1,k_2,x_2}
+(k_1+1)h_{k_1+1,k_2}-(k_2+1)h_{k_1,k_2+1}=0
\]
for any $k_1,k_2$.
Since $h_{k_1,k_2} = 0$ for $k_1 \gg 0$ or $k_2 \gg 0$,
we prove by the descending induction on $(k_1,k_2)$ that 
there exist $h_k \in \mathcal{O}_{X,0}$ for $k \in \mathbb{Z}_{\ge 0}$ 
such that $h_{k_1,k_2} = \binom{k_1+k_2}{k_1}h_{k_1+k_2}$,
and if we expand $h_k$ in a power series
$h_k = \sum_{l_1,l_2} a_{k,l_1,l_2}x_1^{l_1}x_2^{l_2}$,
then we have $a_{k,l_1,l_2} = 0$ unless $l_1=l_2$.
Therefore, $h \in \mathcal{O}_{{Y^{\#}},0}$.

On the other hand, we can solve the equations 
\[
\begin{split}
h_{k_1,k_2} = &x_1g_{k_1,k_2,x_1} - x_2g_{k_1,k_2,x_2} \\
&+(k_1+1)g_{k_1+1,k_2}-(k_2+1)g_{k_1,k_2+1}
\end{split}
\]
by the descending induction on $(k_1,k_2)$
to find 
\[
g = \sum_{k_1,k_2} g_{k_1,k_2}(\log x_1)^{k_1}(\log x_2)^{k_2}
\]
such that $hdx_1/x_1 = dg$ in $\Omega^1_{X^{\#}/Y^{\#}}$.

(2) We check the case where $(X, B) = (\Delta^2, \text{div}(x_1x_2))$, 
$(Y, C) = (\Delta, \{0\})$, $f^*(y) = x_1x_2$ and $p = 0$, 
that is
\[
(\rho_X^1)_*\mathcal{O}_{X^{\#}}
\cong (\text{pr}_1)^{-1}(\mathcal{O}_X) 
\otimes_{(\text{pr}_2)^{-1}\rho_Y^{-1}(\mathcal{O}_Y)}
(\text{pr}_2)^{-1}(\mathcal{O}_{Y^{\#}})
\]
and $R^1(\rho_X^1)_*\mathcal{O}_{X^{\#}} = 0$.
The general case is similar.

We have $\rho_X^{-1}(0) \cong (S^1)^2$ and 
$\rho_Y^{-1}(0) \cong S^1$.
The map $\rho_X^1: \rho_X^{-1}(0) \to \rho_Y^{-1}(0)$
is given by $\rho_X^1(\theta_1, \theta_2) = \theta_1+\theta_2$.

Let $h = \sum_{k_1,k_2} h_{k_1,k_2}(\log x_1)^{k_1}(\log x_2)^{k_2}$ 
be a local section of $\mathcal{O}_{X^{\#}}$, where 
$h_{k_1,k_2} \in \mathcal{O}_{X,0}$.
The monodromies of the multi-valued functions $\log x_1$ and $\log x_2$
along the fibers of $\rho_X^1$ which are homeomorphic to $S^1$'s are given by 
\[
\log x_1 \mapsto \log x_1 + 2\pi\sqrt{-1},
\quad \log x_2 \mapsto \log x_2 - 2\pi\sqrt{-1}.
\]
The function $h$ is single valued along these $S^1$'s if and only if 
there exist $h_k \in \mathcal{O}_{X,0}$ for $k \in \mathbb{Z}_{\ge 0}$ 
such that $h_{k_1,k_2} = \binom{k_1+k_2}{k_1}h_{k_1+k_2}$.
Therefore we have the first equality.

We have an exact sequence
\[
\begin{split}
&0 \to \rho_X^{-1}(\mathcal{O}_X) 
\otimes_{(f^{\#})^{-1}\rho_Y^{-1}(\mathcal{O}_Y)}
(f^{\#})^{-1}(\mathcal{O}_{Y^{\#}}) \\
&\to \mathcal{O}_{X^{\#}} \to
(L_X/(f^{\#})^{-1}L_Y) \otimes \mathcal{O}_{X,0} \to 0
\end{split}
\]
where $L_X$ and $L_Y$ are the locally constant sheaves on 
$\rho_X^{-1}(0) \cong (S^1)^2$ and $\rho_Y^{-1}(0) \cong S^1$
such that the stalks
$L_{X,(0,0)}$ and $L_{Y,0}$ at $(0,0) \in (S^1)^2$ and $0 \in S^1$
have $\mathbb{C}$-bases $\{e_{k_1,k_2}\}_{k_1,k_2 \in \mathbb{Z}_{\ge 0}}$ and
$\{e_k\}_{k \in \mathbb{Z}_{\ge 0}}$
whose the monodromies $M_1, M_2$ around the first and second factors and $M$ 
are given by 
\[
\begin{split}
&M_1(e_{k_1,k_2}) = \sum_{j=0}^{k_1} \binom{k_1}{j}(2\pi\sqrt{-1})^{k_1-j}
e_{j,k_2} \\
&M_2(e_{k_1,k_2}) = \sum_{j=0}^{k_2} \binom{k_2}{j}(2\pi\sqrt{-1})^{k_2-j}
e_{k_1,j} \\ 
&M(e_k) = \sum_{j=0}^{k} \binom{k}{j}(2\pi\sqrt{-1})^{k-j} e_j
\end{split}
\]
and the homomorphism $(f^{\#})^{-1}L_Y \to L_X$ is given by 
\[
(f^{\#})^{-1}e_{k} \mapsto \sum_{j=0}^k \binom{k}{j} e_{j,k-j}.
\]
In other words, we have
\[
\begin{split}
&(M_1 - M_2)(e_{k_1,k_2}) = \sum_{j_1=0}^{k_1} \sum_{j_2=0}^{k_2}
(- 1)^{k_2-j_2} \binom{k_1}{j_1} \binom{k_2}{j_2}
(2\pi\sqrt{-1})^{k_1+k_2-j_1-j_2} e_{j_1,j_2} \\
&(M_1-M_2)((f^{\#})^{-1}e_k) = (f^{\#})^{-1}e_k.
\end{split}
\]
We have
\[
\begin{split}
&(\rho_X^1)_*(f^{\#})^{-1}L_Y \cong (\rho_X^1)_*L_X \cong L_Y \\
&R^1(\rho_X^1)_*(f^{\#})^{-1}L_Y \cong L_Y, \quad
R^1(\rho_X^1)_*L_X = 0
\end{split}
\]
hence 
\[
(\rho_X^1)_*(L_X/f^{\#})^{-1}L_Y) \cong L_Y,
\quad R^1(\rho_X^1)_*(L_X/f^{\#})^{-1}L_Y) \cong 0.
\]
Since 
\[
\begin{split}
&(\rho_X^1)_*\rho_X^{-1}(\mathcal{O}_X) \cong
(\text{pr}_1)^{-1}(\mathcal{O}_X) \\
&R^1(\rho_X^1)_*\rho_X^{-1}(\mathcal{O}_X) \cong 
\mathbb{C}_{(\rho_Y)^{-1}(0)} \otimes \mathcal{O}_{X,0} \\
&\mathcal{O}_{Y^{\#}} \vert_{\rho_Y^{-1}(0)}
\cong L_Y \otimes \mathcal{O}_{Y,0}
\end{split}
\]
we have
\[
\begin{split}
&(\rho_X^1)_*(\rho_X^{-1}(\mathcal{O}_X) 
\otimes_{(f^{\#})^{-1}\rho_Y^{-1}(\mathcal{O}_Y)}
(f^{\#})^{-1}(\mathcal{O}_{Y^{\#}})) \\
&\cong (\text{pr}_1)^{-1}(\mathcal{O}_X)
\otimes_{(\text{pr}_2)^{-1}\rho_Y^{-1}(\mathcal{O}_Y)}
(\text{pr}_2)^{-1}(\mathcal{O}_{Y^{\#}}) \\
&R^1(\rho_X^1)_*(\rho_X^{-1}(\mathcal{O}_X) 
\otimes_{(f^{\#})^{-1}\rho_Y^{-1}(\mathcal{O}_Y)}
(f^{\#})^{-1}(\mathcal{O}_{Y^{\#}})) \\
&\cong L_Y \otimes \mathcal{O}_{X,0}.
\end{split}
\]
Therefore, we have the desired result.

(3) is obtained by combining (1) and (2).
\end{proof}

We assume that $f$ is weakly semistable in the rest of this section, 
and use the notation of \S 3.
Let $m'$ be the number of irreducible components of $C$ which contain $y \in Y$,and let $y_1, \dots, y_{m'}$ be the corresponding local coordinates.
We fix $\bar y \in \rho_Y^{-1}(y)$ and 
denote by $\rho_{\bar D}$ and $\rho_{\bar D_{i_0, \dots, i_t}}$
the restrictions of $\rho_X$ to $\bar D$ and $\bar D_{i_0, \dots, i_t}$.

\begin{Defn}
We define the structure sheaves of $\bar D$ and 
$\bar D_{i_0, \dots, i_t}$ by the following formula:
\[
\begin{split}
&\mathcal{O}_{\bar D} = (\mathcal{O}_{X^{\#}} 
\otimes_{\rho_X^{-1}(\mathcal{O}_X)} \rho_X^{-1}(\mathcal{O}_E))/I \\
&\mathcal{O}_{\bar D_{i_0, \dots, i_t}} = (\mathcal{O}_{X^{\#}} 
\otimes_{\rho_X^{-1}(\mathcal{O}_X)} \rho_X^{-1}
(\mathcal{O}_{E_{i_0, \dots, i_t}}))/I_{i_0, \dots, i_t}
\end{split}
\]
where $I$ and $I_{i_0, \dots, i_t}$ are ideals 
generated by $\log y_j - c_i$ ($j = 1, \dots, m'$) for some constants
$c_i \in \mathbb{C}$ corresponding to the choice of 
$\bar y \in \rho_Y^{-1}(y)$.
We define
\[
\Omega^p_{E/\mathbb{C}}(\log) = 
\Omega^p_{X/Y}(\log) \otimes_{\mathcal{O}_X} \mathcal{O}_E
\]
and define the {\em log de Rham complex} 
on the fiber by 
\[
\begin{split}
&\Omega^p_{\bar D} =
\rho_{\bar D}^{-1}(\Omega^p_{E/\mathbb{C}}(\log)) 
\otimes_{\rho_{\bar D}^{-1}(\mathcal{O}_E)} \mathcal{O}_{\bar D} 
= \Omega^p_{X^{\#}/Y^{\#}} \otimes_{\mathcal{O}_{X^{\#}}} 
\mathcal{O}_{\bar D} \\
&\Omega^p_{\bar D_{i_0, \dots, i_t}} =
\rho_{\bar D}^{-1}(\Omega^p_{E/\mathbb{C}}(\log)) 
\otimes_{\rho_{\bar D}^{-1}(\mathcal{O}_E)} 
\mathcal{O}_{\bar D_{i_0, \dots, i_t}}
= \Omega^p_{X^{\#}/Y^{\#}} \otimes_{\mathcal{O}_{X^{\#}}} 
\mathcal{O}_{\bar D_{i_0, \dots, i_t}}
\end{split}
\]
for $p \ge 0$.
\end{Defn}

We note that there is no ideal sheaf of $D$
corresponding to $I$ because of the monodromies.
The following is easy:

\begin{Lem}\label{Omega^1_{E/C}(log)}
Let $t' = \dim E - \dim E_{i_0, \dots, i_t}$ as in 
Proposition~\ref{direct product}.

(1) If $t=t'=0$, then 
\[
\Omega^p_{E/\mathbb{C}}(\log) \otimes_{\mathcal{O}_E}
\mathcal{O}_{E_{i_0}} \cong \Omega^p_{E_{i_0}}(\log G_{i_0}).
\]

(2) If $t' > 0$, then
\[
\begin{split}
&0 \to \Omega^1_{E_{i_0, \dots, i_t}}(\log G_{i_0, \dots, i_t}) 
\to \Omega^1_{E/\mathbb{C}}(\log)
\otimes_{\mathcal{O}_{E}} \mathcal{O}_{E_{i_0, \dots, i_t}} \\
&\to \mathcal{O}_{E_{i_0, \dots, i_t}}^{\oplus t'} 
\to 0
\end{split}
\]
where the arrow next to the last is the residue homomorphism.
\end{Lem}

The following is similar to Theorem~\ref{relative log Poincare}:

\begin{Thm}\label{log Poincare on the fiber}
(1) 
\[
\begin{split}
&\mathbb{C}_{\bar D} \cong \Omega^{\bullet}_{\bar D} \\
&\mathbb{C}_{\bar D_{i_0, \dots, i_t}} 
\cong \Omega^{\bullet}_{\bar D_{i_0, \dots, i_t}}
\end{split}
\]
in the derived category of sheaves on $\bar D$.

(2)
\[
\begin{split}
&R(\rho_{\bar D})_*\Omega^p_{\bar D}
\cong \Omega^p_{E/\mathbb{C}}(\log) \\
&R(\rho_{\bar D_{i_0, \dots, i_t}})_*\Omega^p_{\bar D_{i_0, \dots, i_t}}
\cong \Omega^p_{E/\mathbb{C}}(\log) \otimes_{\mathcal{O}_E} 
\mathcal{O}_{E_{i_0, \dots, i_t}}
\end{split}
\]
in the derived category of sheaves on $E$.

(3) 
\[
\begin{split}
&R(\rho_{\bar D})_*\mathbb{C}_{\bar D}
\cong \Omega^{\bullet}_{E/\mathbb{C}}(\log) \\
&R(\rho_{\bar D_{i_0, \dots, i_t}})_*\mathbb{C}_{\bar D_{i_0, \dots, i_t}} 
\cong \Omega^{\bullet}_{E/\mathbb{C}}(\log) \otimes_{\mathcal{O}_E} 
\mathcal{O}_{E_{i_0, \dots, i_t}}
\end{split}
\]
in the derived category of sheaves on $E$.
\end{Thm}

%%%%%%%%%%%%%%%%%%%%%%%%%%%%%%%%%%%%%%%%%%%%%%%%%%%%%%%%%%%%%%%%%%

\section{Cohomological mixed Hodge complex}

We prove the main result (Theorem~\ref{main}) of this paper.
We assume that $f$ is weakly semistable in this section, and use the 
notation of \S3.
We shall construct a cohomological mixed Hodge complex on $E = f^{-1}(y)$.

On the $\mathbb{Z}$-level, we put a {\it weight filtration} $W_q$ on 
the complex of sheaves $R(\rho_X)_*\mathbb{Z}_{\bar D}$
on $E$ as in Definition~\ref{weight filtration} by
\[
\begin{split}
&W_q(R(\rho_X)_*\mathbb{Z}_{\bar D}) \\
&= (\tau_{\le q}(R(\rho_X)_*\mathbb{Z}_{\bar D^{[0]}}) \to 
\tau_{\le q+1}(R(\rho_X)_*\mathbb{Z}_{\bar D^{[1]}}) \to \\
&\quad \dots \to \tau_{\le q+t}(R(\rho_X)_*\mathbb{Z}_{\bar D^{[t]}}) 
\to \dots)
\end{split}
\]
in the derived category of complexes of sheaves on 
$E$.

We have an exact sequence
\[
0 \to \mathcal{O}_E \to \mathcal{O}_{E^{[0]}} \to \mathcal{O}_{E^{[1]}}
\to \dots \to \mathcal{O}_{E^{[t]}} \to \dots 
\]
as before, where we identified the sheaves $\mathcal{O}_{E^{[t]}}$
with their direct image sheaves on $E$.

On the $\mathbb{C}$-level, 
we define the weight filtration and the Hodge filtration on the complex
$\Omega^{\bullet}_{E/\mathbb{C}}(\log) = 
\Omega^{\bullet}_{X/Y}(\log) \otimes_{\mathcal{O}_X} \mathcal{O}_E$.

\begin{Defn}
The {\it weight filtration} $\{W_q\}$ on 
$\Omega^{\bullet}_{E/\mathbb{C}}(\log)$ 
is defined by
\[
\begin{split}
&W_q(\Omega^{\bullet}_{E/\mathbb{C}}(\log)) \\
&= (W_q(\Omega^{\bullet}_{E/\mathbb{C}}(\log) 
\otimes \mathcal{O}_{E^{[0]}}) 
\to W_{q+1}(\Omega^{\bullet}_{E/\mathbb{C}}(\log)
\otimes \mathcal{O}_{E^{[1]}}) \to \\
&\quad \dots \to W_{q+t}(\Omega^{\bullet}_{E/\mathbb{C}}(\log)
\otimes \mathcal{O}_{E^{[t]}}) 
\to \dots),
\end{split}
\]
where $W$ on the right hand side is the filtration defined by the order of
log poles at \cite{Deligne}~3.1.5.

The {\it Hodge filtration} $\{F^p\}$ on $\Omega^{\bullet}_{X/Y}(\log)$
is defined by
\[
F^p(\Omega^{\bullet}_{X/Y}(\log)) = 
\sigma_{\ge p}(\Omega^{\bullet}_{X/Y}(\log))
\]
where $\sigma_{\ge p}$ denotes the stupid filtration 
(cf. 1.4.7 of \cite{Deligne}). 
We consider its restriction $\{F^p\}$ on the fiber $E$:
\[
F^p(\Omega^{\bullet}_{E/\mathbb{C}}(\log)) =
\sigma_{\ge p}(\Omega^{\bullet}_{E/\mathbb{C}}(\log)).
\]
\end{Defn}

The following is the main result of this paper:

\begin{Thm}\label{main}
\[
\{(R(\rho_X)_*\mathbb{C}_{\bar D}, W), 
(\Omega^{\bullet}_{E/\mathbb{C}}(\log), W, F)\}
\]
is a cohomological mixed Hodge complex on $E$ (\cite{Deligne}~8.1.6).
\end{Thm}

The proof is reduced to Propositions~\ref{main1} and \ref{main2}. 

\begin{Prop}\label{main1}
\[
(R(\rho_X)_*\mathbb{C}_{\bar D}, W)
\cong (\Omega^{\bullet}_{E/\mathbb{C}}(\log), W)
\]
in the derived category of filtered complexes of 
sheaves on $E$.
\end{Prop}

\begin{proof}
Since the truncation is a canonical functor, this is a consequence of 
Theorem~\ref{log Poincare on the fiber}~(3) and \cite{Deligne}~3.1.8.
\end{proof}

The following is similar to Proposition~\ref{Gr}:

\begin{Lem}
\[
\begin{split}
&\text{Gr}_q^W(\Omega^{\bullet}_{E/\mathbb{C}}(\log)) \\
&\cong \text{Gr}_q^W(\Omega^{\bullet}_{E/\mathbb{C}}(\log) 
\otimes \mathcal{O}_{E^{[0]}})
\oplus \text{Gr}_{q+1}^W(\Omega^{\bullet}_{E/\mathbb{C}}(\log) 
\otimes \mathcal{O}_{E^{[1]}})[-1] \oplus \\
&\quad \dots \oplus \text{Gr}_{q+t}^W(\Omega^{\bullet}_{E/\mathbb{C}}(\log)
\otimes \mathcal{O}_{E^{[t]}})[-t]
\oplus \dots 
\end{split}
\]
\end{Lem}

\begin{Lem}\label{Poincare residue}
Let $(X, B)$ and $(E, G)$ be as in 
Proposition~\ref{real oriented blowing-up of a stratum}.
Then the Poincar\'e residue induces the following isomorphisms
of filtered complexes:
\[
\begin{split}
&(\text{Gr}_q^W(\Omega^{\bullet}_X(\log B)), F)
\cong (\Omega^{\bullet}_{B^{[q]}}[- q], F[-q]) \\
&= (\bigoplus_{1 \le i_1 < \dots < i_q \le N}  
\Omega^{\bullet}_{B_{i_1, \dots, i_q}}[- q], F[-q]) \\
&(\text{Gr}_q^W(\Omega^{\bullet}_X(\log B) 
\otimes_{\mathcal{O}_X} \mathcal{O}_E), F) \cong 
(\bigoplus_{t=0}^q(\Omega^{\bullet}_{G^{[t]}})^{\binom l{q-t}}[- q], F[-q]).
\end{split}
\]
\end{Lem}

We note that $G^{[0]} = E$ by definition.

\begin{proof}
The first formula is \cite{Deligne}~3.1.5.2.
The second is similar.
\end{proof}

\begin{Prop}\label{main2}
\[
(\text{Gr}_q^W(\Omega^{\bullet}_{E/\mathbb{C}}(\log)), F)
\]
is a cohomological Hodge complex of weight $q$ on $E$ (\cite{Deligne}~8.1.2).
\end{Prop}

\begin{proof}
We have to prove that
\[
\begin{split}
&(H^n(E, \text{Gr}_{q+t}^W(\Omega^{\bullet}_{E/\mathbb{C}}(\log)
\otimes \mathcal{O}_{E^{[t]}})[-t]), F) \\
&\cong (H^{n-t}(E, \text{Gr}_{q+t}^W(\Omega^{\bullet}_{E/\mathbb{C}}(\log)
\otimes \mathcal{O}_{E^{[t]}})), F)
\end{split}
\]
is a Hodge structure of weight $n+q$: 
\[
H^{n-t}(E, \text{Gr}_{q+t}^W(\Omega^{\bullet}_{E/\mathbb{C}}(\log)
\otimes \mathcal{O}_{E^{[t]}})) = \sum_p F^p \cap \bar F^{n+q-p}.
\]
By Lemmas~\ref{Omega^1_{E/C}(log)} and \ref{Poincare residue}, we can write
\[
(\text{Gr}_{q+t}^W(\Omega^{\bullet}_{E/\mathbb{C}}(\log)
\otimes \mathcal{O}_{E^{[t]}}), F)
\cong (\sum_G \Omega^{\bullet}_G[-q-t], F[-q-t])
\]
for some set of the strata $G$.
Thus the problem is reduced to prove
\[
H^{n-q-2t}(E, \Omega^{\bullet}_G) = 
\sum_p F^{p-q-t} \cap \bar F^{n-p-t}.
\]
But this is the usual Hodge theorem on $G$.
\end{proof}

By \cite{Deligne}~Scholie 8.1.9, we obtain the following corollary 
(cf. \cite{Fujisawa}):

\begin{Cor}
(1) The spectral sequence
\[
_WE_1^{p,q} = H^{p+q}(E, \text{Gr}_{-p}^W(R(\rho_X)_*\mathbb{C}_{\bar D}))
\Rightarrow H^{p+q}(\bar D, \mathbb{C})
\]
degenerates at $E_2$.

(2) The spectral sequence
\[
_FE_1^{p,q} = H^q(E, \Omega^p_{E/\mathbb{C}}(\log))
\Rightarrow H^{p+q}(\bar D, \mathbb{C})
\]
degenerates at $E_1$.
\end{Cor}

Combining with Corollary~\ref{locally constant sheaf}, 
we obtain the following (\cite{Kollar} and \cite{Nakayama}):

\begin{Cor}
$R^qf_*\Omega^p_{X/Y}(\log)$ is a locally free sheaf for any $p,q$.
\end{Cor}

%%%%%%%%%%%%%%%%%%%%%%%%%%%%%%%%%%%%%%%%%%%%%%%%%%%%%%%%%%%%%%%%%%%%%%%%%%

Department of Mathematical Sciences, University of Tokyo, 

Komaba, Meguro, Tokyo, 153-8914, Japan 

kawamata@ms.u-tokyo.ac.jp

\end{document}